\documentclass[12pt, final]{amsproc}
\usepackage{amssymb,amsbsy,amsmath,amsfonts,mathpazo}
\usepackage{latexsym,euscript,exscale}
\usepackage{mathrsfs}
\usepackage{url}

\usepackage{graphicx}
\usepackage{wrapfig}

\usepackage{amscd}
\usepackage[all]{xy}
\usepackage[notcite,notref]{showkeys}
\theoremstyle{plain}
\newtheorem{theorem}{Theorem}
\newtheorem{proposition}[theorem]{Proposition}
\newtheorem{lemma}[theorem]{Lemma}
\newtheorem{fact}{Fact}
\newtheorem{corollary}[theorem]{Corollary}

\theoremstyle{remark}
\newtheorem{remark}[theorem]{Remark}

\newtheorem{example}[theorem]{Example}
\newtheorem{definition}[theorem]{Definition}

\def\dist{\operatorname{dist}}
\def\dens{\operatorname{dens}}

\def\Lip{\operatorname{Lip}}
\def\N{\mathbb{N}}
\def\R{\mathbb{R}}

\def\U{\mathscr{U}}
\def\K{\mathscr{K}\!}

\def\e{\varepsilon}

\title{Accessible operators on ultraproducts of Banach spaces}
\author{F\'elix Cabello S\'anchez}
\address{Departamento de Matem\'{a}ticas, UEx, 06006-Badajoz, Spain}
\email{fcabello@unex.es}

\thanks{Partially supported by PID2019-103961GB-C21 and Junta de Extremadura IB16056.}

\begin{document}

\begin{abstract}
We address a question by Henry Towsner about the possibility of representing linear operators between ultraproducts of Banach spaces by means of ultraproducts of nonlinear maps.

We provide a bridge between these ``accessible'' operators and the theory of twisted sums through the so-called quasilinear maps. Thus, for many pairs of Banach spaces $X$ and $Y$, there is an ``accessible'' operator $X_\U\to Y_\U$ that is not the ultraproduct of a family of operators $X\to Y$ if and only if there is  a  short exact sequence of quasi-Banach spaces and operators $0\to Y\to Z\to X\to 0$ that does not split.

We then adapt classical work by Ribe and Kalton--Peck to exhibit pretty concrete examples of accessible functionals and endomorphisms for the sequence spaces $\ell_p$.

The paper is organized so that the main ideas are accessible 
to readers working on ultraproducts and requires only a rustic knowledge of Banach space theory.
\end{abstract}

\noindent{{\footnotesize The final version will appear in Extracta Mathematic\ae, \url{https://revista-em.unex.es/index.php/EM} }}\\[10pt]

\maketitle


\section{Introduction}

\subsection{Purpose}
This paper stems from a question posted on Mathoverflow by Henry Towsner and the subsequent comments by the usual suspects. The thread was about the possibility of representing linear functionals on ultraproducts of Banach spaces by means of ultraproducts of nonlinear functions, see \cite{T}.

For reasons that will become clear later we will take a slightly more general approach: we consider families of nonlinear mappings $f_i:X_i\to Y_i$ acting between (quasi-) Banach spaces and we explore the possibility that they define an ``accessible'' linear operator $f:(X_i)_\U\to (Y_i)_\U$ through the formula
\begin{equation}\label{eq:fixi}
f((x_i)_\U)=(f_i(x_i))_\U,
\end{equation}
where $\U$ is a suitable ultrafilter on the underlying index set.
Of course this has little interest if $f$ is an ultraproduct of linear operators.
It turns out that the linear character of \eqref{eq:fixi} depends on the asymptotic behavior of the ``quasilinearity constants''
$$
Q[f_i]=\sup_{x,y\in X_i}\frac{\|f_i(x+y)-f_i(x)-f_i(y)\|}{\|x\|+\|y\|},
$$
which is a very fortunate coincidence because quasilinear maps are a classical tool in the theory of twisted sums of (quasi-) Banach spaces and have been extensively studied in connection with the ``three-space'' problem.

Thus, as we shall see, for many pairs of (quasi-) Banach spaces $X$ and $Y$, one has an ``accessible'' operator $X_\U\to Y_\U$ that is not an ultraproduct of operators if and only if there is a short exact sequence of (quasi-) Banach spaces and operators
$$
\begin{CD}
0@>>> Y @>>> Z @>>> X@>>> 0
\end{CD}
$$
that does not split. This is the leading idea of the paper, which we will exploit in both directions.

Let us explain the organization of the paper and highlight some results. This Section contains, apart from this general introduction, some general definitions and conventions which are used along the paper.

In Section~\ref{main} we characterize those families that define accessible operators by means of their quasilinearity constants, we relate quasilinear maps and accessible operators and we include the minimal background on the ``three-space'' problem that is required to understand the main examples.  

In Section~\ref{sec:appl} we use two old examples of quasilinear maps, due to Ribe and Kalton--Peck, to exhibit pretty explicit examples of accessible functionals  on ultrapowers of $\ell_1$ and accessible endomorphisms on $(\ell_p)_\U$ for each $0<p<\infty$.

Then we present a general procedure which allows one to construct an accessible operator starting with a single quasilinear map.

We also relate Towsner's ideas to the classical notion of a $\mathcal K$-space and we investigate these a bit further. 


Those readers that are primarily interested in concrete examples of accessible operators on ultraproducts of concrete Banach spaces (the sequence spaces $\ell_p$, in fact) can skip the remainder of this introduction, read Sections 2.1 and 2.2 and then go to Sections 3.1--3.3. These sections have been written with all the calculations necessary to understand the examples presented, even if they may appear unnecessary and\,/\,or repetitive to the average Banach spacer who may be familiar with the papers \cite{ribe, kp}.

Although these notes do not contain any essentially new idea, I believe that the interpretation of quasilinear maps as ``accessible operators on ultraproducts'' it provides has some merit.

\subsection{Preliminaires}
We now gather some general definitions and conventions that we use throughout the paper, mainly to fix the notations. The ground field is $\R$, the real numbers, unless otherwise stated.

\subsection*{Quasi-Banach spaces} General references are \cite{kpr, khandbook}.  A quasinorm on a linear space $X$ is a mapping $x\in X\mapsto\|x\|\in \R$ satisfying the following conditions:
\begin{itemize}
\item $\|x\|\geq 0$ for every $x\in X$ and $\|x\|=0$ if and only if $x=0$.
\item $\|tx\|=|t|\|x\|$ for every $t\in\mathbb R$ and every $x\in X$.
\item There is a constant $\Delta$ such that $\|x+y\|\leq\Delta(\|x\|+\|y\|)$ for every $x,y\in X$.
\end{itemize}
The least possible constant $\Delta$ fitting in the preceding inequality is often called the ``modulus of concavity'' of the quasinorm. A quasinormed space is a linear space equipped with a quasinorm. On such a space $X$ there is a linear topology for which the ``balls''
$$
B(r)=\{x\in X: \|x\|\leq r\} 
$$
form a neighborhood base at zero. If $X$ is complete we call it a quasi-Banach space. It is important to realize that a quasinorm may be discontinuous for the topology it induces! To avoid these complications we will consider mostly $p$-norms ($0<p\leq 1$), that is, quasinorms satisfying the inequality $\|x+y\|^p\leq \|x\|^p+\|y\|^p$ for every $x,y\in X$. The modulus of concavity of a $p$-norm is at most $2^{1/p-1}$. The Aoki-Rolewicz theorem states that every quasinorm is equivalent to a $p$-norm for some $0<p\leq 1$, where $p$ can be taken so that $2=(2\Delta)^p$.
Note that if $0<q<p\leq 1$, then every $p$-norm is also a $q$-norm.
\smallskip

The simplest examples of quasi-Banach spaces are the sequence spaces $\ell_p$ for $0<p<\infty$. 
In general we regard (numerical) sequences as functions $x:\N\to\mathbb R$, so that $x$ belongs to $\ell_p$ if the quasinorm
$$
\|x\|_p=\left( \sum_{n=1}^\infty |x(n)|^p\right)^{1/p}
$$
is finite. For $p\geq 1$, the function $\|\cdot\|_p$ is a norm, but it is just a $p$-norm when $0<p<1$. We denote by $e_k$ the sequence that takes the value 1 at $k$ and 0 elsewhere.

\smallskip

Other quasi-Banach spaces that will play a secondary role in these notes are the Lebesgue spaces $L_p$, the Hardy classes $H_p$ (see \cite{kpr}), the $p$-Gurariy space $G_p$ (see \cite{cgk} or \cite[Chapter~6]{HMBST}) and the noncommutative $L^p$-spaces \cite{p-x}.

\subsection*{Ultraproducts} General references are \cite{heinrich, sims, k84}.  
We consider ultraproducts of quasi-Banach spaces ``inside'' the set-theoretic ultraproduct --- which will be needed in Section~\ref{sec:AFNLCS}. So, let $X_i$ be a family of sets indexed by $I$ and let $\U$ be a countably incomplete ultrafilter on $I$. The set-theoretic ultraproduct $\langle X_i\rangle_\U$ is the product $\prod_{i\in I} X_i$ factored by the equivalence relation that identifies $(x_i)$ and $(y_i)$ if the set $\{i\in I: x_i=y_i\}$ belongs to $\U$. The class of $(x_i)$ in $\langle X_i\rangle_\U$ shall be denoted by $\langle x_i\rangle_\U$. 
Note that each family of maps $f_i:X_i\to Y_i$ induces a map 
$\langle f_i\rangle_\U : \langle X_i\rangle_\U\to \langle Y_i\rangle_\U$ sending $\langle x_i\rangle_\U$
to $\langle f_i(x_i)\rangle_\U$.

When each $X_i$ is a quasinormed space with quasinorm $\|\cdot\|_i$ one can measure the ``size'' of $\langle x_i\rangle_\U$ using the ``limit functional'' $ \langle X_i\rangle_\U\to [0, \infty]$ given by
\begin{equation}\label{quasiU}
\langle x_i\rangle_\U\longmapsto \lim_{\U(i)}\|x_i\|
\end{equation}
Then one naturally wants to consider the subset of ``finite'' elements of $\langle X_i\rangle_\U$ modulo ``infinitesimals'' to obtain the ``quasi-Banach'' ultraproduct, but a serious problem appears: the sum of two ``finite'' elements need not to be ``finite'' if the moduli of concavity are not uniformly bounded.  The simplest way to avoid these pathologies is to consider only families of $p$-Banach spaces for some fixed $0<p\leq 1$. Indeed, if each $X_i$ is $p$-normed, then 
$$
\langle X_i\rangle_\U^{\text{fin}}=\left\{ \langle x_i\rangle_\U: \lim_{\U(i)}\|x_i\|_i<\infty  \right\}\quad\quad\text{and}\quad\quad
\langle X_i\rangle_\U^{\text{inf}}=\left\{ \langle x_i\rangle_\U: \lim_{\U(i)}\|x_i\|_i=0  \right\}
$$
are linear subspaces of  $\langle X_i\rangle_\U$ and the functional
(\ref{quasiU}) defines a quasinorm (a $p$-norm in fact) on the quotient
$
\langle X_i\rangle_\U^{\text{fin}}/\langle X_i\rangle_\U^{\text{inf}}
$.
The resulting $p$-normed space (which is automatically complete) is the quasi-Banach space ultraproduct of the family $X_i$ with respect to $\U$ and shall be denoted by $(X_i)_\U$.

Of course this is equivalent to consider the quotient of the space of bounded families
$$
\ell_\infty(I, X_i)=\left\{(x_i)\in\prod_i X_i: \sup_i\|x_i\|_i<\infty\right\}
$$
(equipped with the obvious $p$-norm) by the subspace
$$c_0^\U(I, X_i)=\{(x_i)\in \ell_\infty(I, X_i): \|x_i\|_i\to 0 \text{ along } \U\}.$$ 
Anyway the fact that the quotient quasinorm on $\ell_\infty(I, X_i)/c_0^\U(I, X_i)$ agrees with the ultraproduct quasinorm in (\ref{quasiU}) depends on the continuity of $p$-norms.

From now on the class of a bounded family $(x_i)$ in $(X_i)_\U$ will be denoted by $(x_i)_\U$. Note that $(x_i)_\U$ ``is'' $\langle x_i\rangle_\U$ ``modulo infinitesimals''.

\subsection*{Homogeneous maps and distances}
 A mapping $f:X\to Y$ acting between linear spaces (over $\R$) is homogeneous if $f(tx)=tf(x)$ for every $t\in\R$ and every $x\in X$.
 
 If $X$ and $Y$ are quasinormed spaces, then $f:X\to Y$ is said to be bounded if it obeys an estimate of the form $\|f(x)\|\leq M\|x\|$ for some constant $M$ and all $x\in X$. A homogeneous map is bounded if and only if the expression
\begin{equation}\label{f}
\|f\|=\sup_{\|x\|=1}\|f(x)\|
\end{equation}
 is finite. This notation is coherent with the usual one for bounded operators. If $f$ and $g$ are homogeneous maps acting between the same quasinormed spaces, then the (possibly infinite) ``distance'' between $f$ and $g$ is defined to be $\|f-g\|$. Also, if $G$ is a family of homogeneous maps, we put
 $$
 \dist(f,G)=\inf_{g\in G}\|f-g\|.
 $$
 Given quasinormed spaces $X$ and $Y$, the space of bounded linear maps $f:X\to Y$ is denoted by $L(X,Y)$. This is a quasinormed space with the quasinorm given in (\ref{f}), which makes it complete when $Y$ is. When $Y=\R$ one obtains the dual $X'$, which is a Banach space even if $X$ is just a quasinormed space.

\section{Accessible operators and quasilinear maps}\label{main}
\subsection{Admissible families and accessible maps}
Let $\U$ be a countably incomplete ultrafilter on $I$ that we consider fixed in all what follows. We refer to $\U$ the notions of ``almost every'' $i\in I$, limits of functions defined on $I$, ``admissible family'', ``accessible operator'', and the like unless otherwise stated.

\begin{definition}
A family of mappings $f_i: X_i\to Y_i$ indexed by $I$ is {\bf admissible} if there is a mapping $f:(X_i)_\U\to (Y_i)_\U$ such that $f((x_i)_\U)=(f_i(x_i))_\U$. The map $f$ itself shall be termed {\bf accessible} and can be properly denoted by $(f_i)_\U$.
\end{definition}

From the perspective of the set-theoretic ultraproduct the admissibility of $(f_i)$ means that $\langle f_i\rangle_\U$ takes finite elements of $\langle X_i\rangle_\U$ into finite elements of $\langle Y_i\rangle_\U$ and that $\langle f_i( x_i)- f_i( y_i)\rangle_\U$ is infinitesimal when $\langle x_i- y_i\rangle_\U$ is.
The latter condition means that $(f_i)$ has the following continuity property: if $\|x_i-y_i\|_i\to0$ along $\U$, then $\|f(x_i)-f(y_i)\|_i\to0$ along $\U$.

Two admissible families $(f_i)$ and $(g_i)$ are said to be equivalent (with respect to $\U$) if $(f_i)_\U= (g_i)_\U$, as mappings from $(X_i)_\U$ to $(Y_i)_\U$.

Since we are interested in families that give rise to operators it will be convenient to consider only uniformly bounded families of homogeneous maps. This will simplify the exposition without any loss of generality, as we now see. Given an arbitrary mapping $u:X\to Y$, define $\tilde u: X\to Y$ by
$$
\tilde u(x)=\frac{\|x\|}{2}\left( u\left( \frac{x}{\|x\|}	\right) - u\left( \frac{-x}{\|x\|}	\right)		 \right)
$$
for $x\neq 0$, and $\tilde u(0)=0$. Obviously $\tilde u$ is homogeneous and, actually, $\tilde u=u$ if and only if $u$ is homogeneous.

\begin{lemma}\label{unif}
Assume that $f_i: X_i\to Y_i$ form an admissible family. If $(f_i)_\U$ is homogeneous, then:
\begin{itemize}
\item $(f_i)$ is equivalent to $(\tilde f_i)$.
\item $(f_i)_\U$ is bounded.
\end{itemize}
Hence every accessible homogeneous map is induced by a uniformly bounded family of homogeneous maps.
\end{lemma}

\begin{proof}
For the first part, just observe that $(f_i)_\U=\widetilde{(f_i)_\U}=(\tilde f_i)_\U$. For the second, consider the bounds $\|\tilde f_i\|$. Clearly, $\|(f_i)_\U\|=  \|(\tilde f_i)_\U\|= \lim_{\U(i)}\|\tilde f_i\|$. This limit has to be finite, for if not, we may choose, for each $i\in I$, a normalized $x_i\in X_i$ such that $\|\tilde f_i(x_i)\|\geq{1\over 2} \|\tilde f_i\|$ and the class of $(f_i(x_i))$ cannot fall in $(Y_i)_\U$.
\end{proof}

\subsection{Quasilinear maps and accessible operators}
We now arrive at the key notion to characterize those families that give rise to linear operators between the corresponding ultraproducts.

\begin{definition}\label{quasi}
A mapping $f:X\to Y$ acting between quasinormed spaces is said to be quasilinear if it is homogeneous and there is a constant $Q$ such that $$\|f(x+y)-f(x)-f(y)\|\leq Q(\|x\|+\|y\|)$$ for all $x,y\in X$.
\end{definition}

We denote by $Q[f]$ the smallest constant $Q$ fitting in the preceding inequality and call it the quasilinearity constant of $f$.

We say that $f$ is trivial if there is a linear map $\ell:X\to Y$ such that $\|f-\ell\|<\infty$. In general, the ``approximating'' map $\ell$ is not bounded. It is important to realize, however, that quasilinear maps defined on finite-dimensional spaces are automatically bounded and that any map whose  
distance to a bounded map is finite is itself bounded. We refer the reader to \cite[Chapter~3]{HMBST} for an account on the basic properties of quasilinear maps and their applications. After all, we are here to sell our book.

\begin{theorem}\label{char}
Let $f_i: X_i\to Y_i$ be a family of homogeneous mappings acting between $p$-Banach spaces indexed by $I$ and let $\U$ be a countably incomplete ultrafilter on $I$. The following are equivalent:
\begin{itemize}
\item[(a)] $(f_i)$ is admissible and $(f_i)_\U:(X_i)_\U \to  (Y_i)_\U$ is a bounded linear operator.

\item[(b)] $\lim_{\U(i)}\|f_i\|<\infty$ and $\lim_{\U(i)}Q[f_i]=0$.
\end{itemize}
Moreover, $(f_i)_\U$ is an ultraproduct of operators if and only if there is a family $u_i\in L(X_i,Y_i)$ such that $\|f_i-u_i\|\to 0$ along $\U$.
\end{theorem}

\begin{proof} Recall our convention that all limits are taken with respect to $\U$. We first prove that (b) implies (a). To show that $(f_i)$ is admissible let us first observe that if $(x_i)$ is a family such that $\|x_i\|\longrightarrow 0$, then  $\|f_i(x_i)\|\longrightarrow 0$.

Now suppose that $(x_i)$ and $(y_i)$ are bounded families such that $\|x_i-y_i\|\to 0$.
Then one has
$$
\|f_i(x_i-y_i)-f_i(x_i)+f_i(y_i)\|\leq Q[f_i](\|x_i\|+\|y_i\|)\longrightarrow 0 
$$
and since $\|f_i(x_i-y_i)\|\longrightarrow 0$ we have $\|f_i(y_i)-f_i(x_i)\|\to 0$ along $\U$ and so $(f_i)$ is admissible.

The linearity of the induced mapping $(f_i)_\U$ is clear since it is obviously homogeneous and, given bounded families $(x_i)$ and $(y_i)$, one has
$
(f_i(x_i+y_i))_\U= (f_i(x_i)+f(y_i))_\U
$
since
$$
\|f_i(x_i+y_i)-f_i(x_i)-f_i(y_i)\|\leq Q[f_i](\|x_i\|+\|y_i\|)\longrightarrow 0
$$
with respect to $\U$.
\smallskip

To prove the converse, suppose that $(f_i)$ is admissible and $(f_i)_\U$ is linear. We know from Lemma~\ref{unif} that $\|f_i\|$ converges to a finite limit and we must show that $Q[f_i] \to 0$.
For each $i\in I$, pick $x_i, y_i\in X_i$ such that $\|x_i\|+\|y_i\|=2$ and
$$
\|f_i(x_i+y_i)-f_i(x_i)-f_i(y_i)\|\geq Q[f_i].
$$
Now, if $(f_i(x_i+y_i))_\U=(f_i(x_i))_\U+(f_i(y_i))_\U$ we see that $Q[f_i] \to 0$.

\smallskip

The ``moreover'' part is obvious once one realizes that two admissible families of homogeneous maps $(f_i)$ and $(g_i)$ are equivalent if and only if $\|f_i-g_i\|\to 0$.
\end{proof}

Thus everything here depends on the relation between the quasilinearity constant $Q[f_i]$ and the distance from $f_i$ to the space of operators $L(X_i,Y_i)$. Let us introduce the following ``approximation constants'', as in \cite[pp. 30--31]{k-o} and \cite[Definition~1]{2c-isr} (see also \cite[pp.~159 and 232]{HMBST}):

\begin{definition}\label{K0}
Given quasi-Banach spaces $X$ and $Y$ we denote by $K_0[X, Y ]$ the
infimum of those constants $C$ for which the following statement is true: if
$f:X\to Y$ is a bounded quasilinear map, then there is a linear map $\ell:X\to Y$
such that $\|f-\ell\|\leq CQ[f]$.
We define $K[X, Y ]$  analogously, just omitting the word ``bounded''. When $Y=\R$ we just write $K_0[X]$ and $K[X]$.

We say that $X$ is a $\mathcal K$-space if $K[X]$ is finite. If $K_0[X]$ is finite, $X$ is said to be a $\mathcal K_0$-space.
\end{definition}

See Section~\ref{sec:AFBS} for examples of $\mathcal K$-spaces.

\begin{corollary}\label{corK0} Let $(X_i)$ and $(Y_i)$ be families of $p$-Banach spaces indexed by $I$ and let $\U$ be a countably incomplete ultrafilter on $I$. The following statements are equivalent:
\begin{itemize}
\item[(a)] Every accessible operator from  $(X_i)$ to $(Y_i)$ is an ultraproduct of operators.

\item[(b)] $\lim_{\U(i)} K_0[X_i, Y_i]<\infty$.
\end{itemize}

In particular, $K_0[X,Y]$ is finite if and only if each accessible operator $X_\U\to Y_\U$ is an ultraproduct of operators, while $X$ is a $\mathcal K_0$-space if and only if the only accessible elements of $(X_\U)'$ are those in $X'_\U$.
\end{corollary}

\subsection{The ``three-space'' problem} Quasilinear maps were introduced in (quasi-) Banach space theory in connection with the ``three-space'' problem. They arise from, give rise to, and describe, short exact sequences of quasi-Banach spaces and operators. Without entering into many details, suppose $\phi: X\to Y$ is a quasilinear map. Then we can equip the product space $Y\oplus X$ with the quasinorm $\|(y,x)\|_\phi=\|y-\phi(x)\|+\|x\|$ (yes, this is a quasinorm!). The resulting quasinormed space is often denoted $Y\oplus_\phi X$ and the diagram
\begin{equation}\label{sex}
\begin{CD}
0@>>> Y @>\imath>> Y\oplus_\phi X @>\pi>> X@>>> 0,
\end{CD}
\end{equation}
where $\imath(y)=(y,0)$ and $\pi(y,x)=0$, is a short exact sequence. Moreover, both $\imath$ and $\pi$ are bounded and relatively open operators which already implies that $Y\oplus_\phi X$ is complete. Besides, $\phi$ is trivial if and only if (\ref{sex}) splits, that is, there is a bounded operator $u:Y\oplus_\phi X\to Y$ such that $u\circ\imath={\bf 1}_Y$ or equivalently a bounded operator $v:X\to Y\oplus_\phi X$ such that $\pi\circ v={\bf 1}_X$. Here (and throughout the paper) ${\bf 1}_S$ denotes the identity on the set $S$.

Thus, quasilinear maps give rise to exact sequences. And conversely: suppose
\begin{equation}\label{sexo}
\begin{CD}
0@>>> Y @>\jmath>> Z @>\varpi>> X@>>> 0
\end{CD}
\end{equation}
is any exact sequence of quasi-Banach spaces and operators (after ``renormalization'' we may think of $Z$ as a quasi-Banach space containing $Y$ as a subspace in such a way that $\jmath$ is the inclusion map, $Z/Y=X$ and $\varpi$ is the quotient map). Let $B:X\to Z$ be a homogeneous bounded section of the quotient map (that is, $\varpi\circ B={\bf 1}_X$ and $L:X\to Z$ a linear (generally unbounded) section. Then the difference $\phi=B-L$ takes values in $Y=\ker\varpi$ and is quasilinear since (it is obviously homogeneous  and)
$$
\|\phi(x+y)-\phi(x)-\phi(y)\|_Y=\|B(x+y)-B(x)-B(y)\|_Z\leq C(\|x\|_X+\|y\|_X).
$$
Now, the ``twisted sum'' $Y\oplus_\phi X$ is isomorphic to $Z$ through the map $u(y,x)=y+L(x)$ and we have a commutative diagram
\begin{equation*}
\xymatrixcolsep{3.5pc}
\xymatrixrowsep{1.25pc}
\xymatrix{
 && Y\oplus_\phi X  \ar[rd]^\pi \ar[dd]^u\\
0\ar[r] & Y \ar[ru]^\imath \ar[rd]_\jmath  &  & X \ar[r] & 0\\
 && Z  \ar[ru]_\varpi
}
\end{equation*}
As the reader can imagine the starting sequence (\ref{sex}) splits if and only if the quasilinear map $\phi$ is trivial.
All this can be seen in \cite[Chapter~3]{HMBST}. 
These ideas have been extensively used in the study of exact sequences of (quasi-) Banach spaces and there is a lot of work around, mostly by Kalton and co-workers.


For the convenience of the reader we now 
present a pageant of ``three space'' problems:

\smallskip

$\bigstar$ Enflo, Lindenstrauss and Pisier \cite[Section~4]{elp} established the existence of nontrivial quasilinear selfmaps on $\ell_2$ thus showing that there are exact sequences $0\to\ell_2\to Z\to\ell_2\to 0$ in which $Z$ is a Banach space not isomorphic to a Hilbert space.

\smallskip

$\bigstar$ Kalton \cite[Section~4]{k78} and, independently, Ribe \cite{ribe} constructed explicit examples of nontrivial quasilinear functionals on $\ell_1$ thus showing that there exist nonlocally convex quasi-Banach spaces whose quotient by a line is locally convex! The quotient is $\ell_1$ in both cases. Hence $\ell_1$ is not a $\K$-space.

\smallskip

$\bigstar$ Then Kalton and Peck \cite{kp} gave examples of nontrivial quasilinear maps on most sequence spaces, in particular, on $\ell_p$ for $0<p<\infty$. These can be regarded as ``vector-valued'' generalizations of Ribe's and shall be discussed in some detail in Section~\ref{explicit}.

\smallskip

$\bigstar$ Concerning the constants appearing in Definition~\ref{K0}, Kalton proved in \cite[Proposition 3.3]{k78} that if $X$ and $Y$ are quasi-Banach spaces such that every quasilinear map from $X$ to $Y$ is trivial, then $K[X,Y]$ is finite and so is $K_0[X,Y]$. 
It it shown in \cite{2c-isr} that in many cases, including those appearing in the hypothesis of the following result (but not always!), $K[X,Y]$ is finite when $K_0[X,Y]$ is; this can be seen also in \cite[Theorem 5.3.13]{HMBST}. One therefore has the following complement to Corollary~\ref{corK0}:

\begin{corollary}
Let $X$ and $Y$ be quasi-Banach spaces, with $Y$ an ultrasummand. Suppose that either $X$ has the BAP or $X$ is a $\mathcal K$-space and $Y$ has the BAP. Then the following statements are equivalent
\begin{itemize}
    \item Every accessible operator from $X_\U$ to $Y_\U$ is an ultraproduct of operators.
    \item $K[X,Y]$ is finite.
    \item $K_0[X,Y]$ is finite.
    \item Every quasilinear map from $X$ to $Y$ is trivial.
    \item Every exact sequence of quasi-Banach spaces $0\to Y\to Z \to X\to 0$ splits. 
\end{itemize}
\end{corollary}

Here, we say that a quasi-Banach space $X$ has the bounded approximation property (BAP, for short) if there is a constant $\lambda\geq 1$ so that, for every finite-dimensional $E\subset X$ and every $\e>0$ there is a finite-rank operator $t\in L(X)$ such that $\|t(x)-x\|\leq\e\|x\|$ for all $x\in E$, with $\|t\|\leq\lambda$.
The term ``ultrasummand'' is defined in Section~\ref{truncation}.

The last item of the preceding corollary requires leaving the world of Banach spaces and considering quasi-Banach spaces as well. Indeed, if $X$ is either $\ell_1$ or $L_1$ and $Y$ is $\ell_p$ or $L_p$ for $1\leq p \leq \infty$, then every exact sequence of Banach spaces $0\to Y\to Z \to X\to 0$ splits. However, not every accessible operator $X_\U\to Y_\U$ is an ultraproduct of operators ... precisely because there are exact sequences $0\to Y\to Z \to X\to 0$ in which $Z$ is not a Banach space!

\begin{wrapfigure}{r}{0.5\textwidth}
  \begin{center}\vspace{-0pt}
    \includegraphics[width=0.5\textwidth]{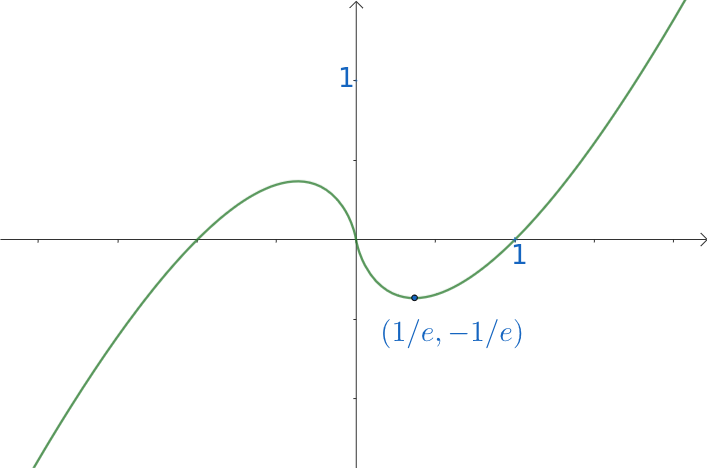}
  \end{center}
  \vspace{-0pt}{\begin{center}{\footnotesize{Graph of Ribe's function $\omega(t)=t\log|t|$.}}\end{center}}
\end{wrapfigure}

\section{Applications}\label{sec:appl}

\subsection{Examples of accessible operators and functionals on $\ell_p$}\label{explicit} To support these ideas, let us consider two concrete examples based on Ribe and Kalton--Peck maps.
The material of this section is taken from \cite{ribe} and \cite{kp} --- and copied, pasted and edited from \cite[Section 3.2]{HMBST}.

 The computations will be explained in detail since we need a very precise control of the distances to the spaces of linear maps to guarantee that the resulting operators are not ultraproducts of linear operators! Most computations depend on the following property of the function $\omega:\R\to \R$ defined by $\omega(t)=t\log|t|$ for $t\neq 0$ and extended by continuity as $\omega(0)=0$. Note that $\omega$ is odd and that $\max\{|\omega(t)|: t\in[0,1]\}=e^{-1}$, attained at $t=e^{-1}$, see the figure just above.

\begin{lemma}[Ribe]\label{w}
For every $s,t\in\R$ one has $|\omega(s+t)-\omega(s)-\omega(t)|\leq (\log 2)(|s|+|t|)$.
\end{lemma}

\begin{proof}
First, suppose $s,t>0$. Then
$$
\frac{|\omega(s+t)-\omega(s)-\omega(t)|}{s+t}=\left|\frac{s}{s+t}\log\frac{s}{s+t} +   \frac{t}{s+t}\log\frac{t}{s+t} \right|= |s'\log s'+ t'\log t'|,
$$
where we have written
 $s'=s/(s+t), t'=t/(s+t)$.
Since $s'+t'=1$  one can use brute-force calculus to obtain that the maximum value of the right-hand side (attained at $s'=t'={1\over 2}$) is $\log 2$, and we are done. If $s,t<0$ just use the odd character of $\omega$ to get the same conclusion.

Now, if $s$ and $t$ have different signs, we may assume that $s$ is positive, $t$ negative and $s+t>0$, so that $|-t|+|s+t|=|s|\leq |s|+|t|$. Since $\omega$ is odd we have
\begin{gather*}
|\omega(s+t)-\omega(s)-\omega(t)|= |\omega(s)-\omega(-t)-\omega(s+t)|
\\
\leq (\log 2)(|-t|+|s+t|)\leq  (\log 2)(|s|+|t|).
\end{gather*}
Nice, isn't it?
\end{proof}

\begin{example}[Ribe's functional]
Let $\ell_1^0$ denote the subspace of finitely supported sequences in $\ell_1$. We define $\varrho:\ell_1^0\to\R$ by the formula
$$
\varrho(x)=\sum_{x(k)\neq 0} x(k)\log\frac{|s(x)|}{|x(k)|},
$$
where $s(x)=\sum_k x(k)$.
Then:\begin{itemize}
\item[(a)] $\varrho$ is quasilinear, with $Q[\varrho]\leq 2\log 2$.
\item[(b)] If we regard $\ell_1^n$ as the subspace of those $x\in\ell_1$ such that $x(k)=0$ for $k>n$, then  ${1\over 2}\log n\leq \dist(\varrho|_{\ell_1^n},(\ell_1^n)')\leq \log n$.
\end{itemize}
Hence, for each free ultrafilter on the integers $\U$, the family $((\log n)^{-1}\varrho)_n$ defines a nontrivial accessible functional on $(\ell_1^n)_\U$. This provides the ``explicit'' formula
$$
f((x_n)_\U)= \lim_{\U(n)}\left\{ \frac{1}{\log n} \sum_{x_n(k)\neq 0} x_n(k)\log\frac{|s(x_n)|}{|x_n(k)|}\right\}.
$$
\end{example}

\begin{proof}(a) Clearly, $\varrho$ is homogeneous.
Note that $\varrho(x)=\omega(s(x))-\sum_k \omega(x(k))$.
Hence, for $x,y\in\ell_1^0$ we have
\begin{align*}
 &|\varrho(x+y)-\varrho(x)-\varrho(y)|\\
=& \left|\omega(s(x+y))-\sum_k \omega(x(k)+y(k))-\omega(sx)+\sum_k \omega(x(k))-\omega(sy)+\sum_k \omega(y(k))\right|   \\
\leq &(\log 2)(|s(x)|+|s(y)|)+ (\log 2)\sum_k \left(|x(k)|+|y(k)|\right)\\
\leq &(2\log 2)(\|x\|_1+\|y\|_1).
\end{align*}
(b) Let $\ell: \ell_1^n\to\R$ be any linear map and set $\delta=\|\varrho|_{\ell_1^n}-\ell\|$. Then since $\varrho(e_i)=0$ for all $i\in\N$ one has $|\ell(e_i)|\leq\delta$. Put $s_n=e_1+\cdots+ e_n$. Then $\|s_n\|_1=n$ and 
$
\left|\varrho(s_n)- \ell\left(s_n\right)  \right|\leq\delta n.
$
But $\varrho(s_n)=n\log n$ and $|\ell(s_n)|=|\sum_i\ell(e_i)|\leq n\delta$, so
$
n\log n\leq 2\delta n
$
and the result follows.
\end{proof}

In all what follows we write $\Lip$ for the space of  Lipschitz functions $\theta: \R\to\R$ such that $\theta(t)=0$ for all $t\leq 0$ and $\|\theta\|_\infty= \sup_t|\theta(t)|<\infty$. The  Lipschitz constant of such an $\theta$ is denoted by $L(\theta)$.

\begin{example}[The Kalton--Peck maps]\label{ex:kpmaps}
Given $\theta\in\Lip$ and $0<p<\infty$ we define $\Phi: \ell_p\to \ell_p$ by the formula
\begin{equation}\label{Fi}
\Phi(x)= x\cdot \theta\left(\log\frac{\|x\|_p}{|x|}\right),
\end{equation}
that is, $\Phi(x)(k)= x(k)\cdot \theta\left(-\log(|x(k)|/{\|x\|_p})\right)$ if $x(k)\neq 0$ and $\Phi(x)(k)=0$ when $x(k)=0$.
Then:
\begin{itemize}
\item[(a)] $\Phi$ is quasilinear on $\ell_p$ and $Q[\Phi]$ depends only on $L(\theta)$ and $p$.
\item[(b)] If, besides, $\theta$ is increasing, then
$
c_p\|\theta\|_\infty\leq \dist(\Phi, L(\ell_p)) \leq \|\Phi:\ell_p\to\ell_p\|\leq \|\theta\|_\infty,
$
where $c_p>0$ is a constant depending only on $p$.
\end{itemize}
Hence:
\begin{itemize}
\item[(c)] If $(\theta_n)$ is a sequence of increasing functions in $\Lip$, with $\sup_n L(\theta_n)<\infty$ and $\lim_n\|\theta_n\|_\infty=\infty$, then, for every free ultrafilter $\U$ on the integers, the formula
    $$
(x_n)_\U\longmapsto \left( \frac{x_n}{\|\theta_n\|_\infty}\cdot \theta_n\left(\log\frac{\|x_n\|_p}{|x_n|} \right) \right)_\U
    $$
defines an accessible endomorphism of $(\ell_p)_\U$ that is not an ultraproduct of operators.
\item[(d)] In particular, letting $\theta_n(t)=\min(t,n)$, we obtain the accessible operator
   $$
(x_n)_\U\longmapsto \left( \frac{x_n}{n}\cdot \theta_n\left(\log\frac{\|x_n\|_p}{|x_n|}\right) \right)_\U=
\left(\left(\left( \frac{x_n(k)}{n}\cdot \min\left\{n, \log\frac{\|x_n\|_p}{|x_n(k)|}\right\}\right)_{\!k\:} \right)_{\!n\:} \right)_\U
    $$
\end{itemize}
\end{example}

\begin{proof}
Recall that for $\omega(t)=t\log|t|$  one has $\sup \{|\omega(t)|: 0\leq t\leq 1\}=e^{-1}$.

(a) We fix $\theta$ and $p$ and we first consider the nonhomogeneous map $\Phi_0:\ell_p\to\ell_p$ given by $\Phi_0(x)=x\cdot\theta(-\log|x|)$, that is:
$$
\Phi_0(x)(k)=\begin{cases} x(k)\cdot\theta(-\log|x(k)|) & \text{if } x(k)\neq 0\\
0 & \text{if } x(k)= 0.\end{cases}
$$
We want to obtain the pointwise estimate
\begin{equation}\label{point}
| \Phi_0(x+y)(k)- \Phi_0(x)(k)- \Phi_0(y)(k)|\leq  2L(\theta)e^{-1}(|x(k)|+|y(k)|)
\end{equation}
To proceed, let us consider the Ribe-like function $\omega_\theta:\R\to \R$ given by $\omega_\theta(t)=t\theta(-\log|t|)$ and prove that
\begin{equation}\label{pointo}
|\omega_\theta(t+s)- \omega_\theta(t)- \omega_\theta(s)|\leq  2L(\theta)e^{-1}(|t|+|s|).
\end{equation}
This implies (\ref{point}) as $\Phi_0(x)=\omega_\theta\circ x$.

If $s,t>0$ we have
\begin{align*}
&\frac{|\omega_\theta(t+s)- \omega_\theta(t)- \omega_\theta(s)|}{ s+t}\\
= & \left| \frac{ (s+t)}{(s+t)}\cdot \theta(-\log|s+t|)-
\frac{ s\cdot  \theta(-\log|s|)}{s+t}
-\frac{t \cdot \theta(-\log|t|)}{s+t}
\right|\\
= &\left|
\frac{ s( \theta(-\log|s+t|)- \theta(-\log|s|))}{s+t}
\right|+ \left|
\frac{ t(\theta(-\log|s+t|)- \theta(-\log|t|))}{s+t}
\right|\\
\leq & 
 L(\theta) \left|
\frac{ s}{s+t}\cdot  \log \frac{s}{s+t}
\right|+ 
 L(\theta) \left|
\frac{ t}{s+t}\cdot  \log \frac{ t}{s+t}
\right| \\
\leq & 2L(\theta)e^{-1}
\end{align*}
The other cases follow from this, as in the proof of Lemma~\ref{w}, which establishes (\ref{pointo}) and (\ref{point}).
Hence,
\begin{align*}
\| \Phi_0(x+y)- \Phi_0(x)- \Phi_0(y)\|_p &\leq  \frac{2L(\theta)}{e}(\||x|+|y|\|_p)
\\
&\leq
\begin{cases}
\dfrac{2L(\theta)}{e}(\|x\|_p+\|y\|_p)& \text{if $p\geq 1$,}\\
\\
\dfrac{2^{1/p}L(\theta)}{e}(\|x\|_p+\|y\|_p) &\text{if $0<p< 1$.}
\end{cases}
\end{align*}
To complete the proof observe that
\begin{equation}\label{eq:estimate}
\|\Phi(x)-\Phi_0(x)\|_p=\biggl{\|}x\cdot\underbrace{\left\{\theta\left(\log\frac{\|x\|_p}{|x|} \right)-  \theta\left(-\log|x|\right)\right\}}_{(\star)}\biggr{\|}_p\leq
L(\theta) \|x\|_p \log\|x\|_p 
\end{equation}
since $(\star)$ is ``pointwise dominated'' by $L(\theta)\log\|x\|_p$.
In particular, for $\|x\|_p\leq 1$ we have 
\begin{equation}\label{p-p0}
\|\Phi(x)-\Phi_0(x)\|_p\leq e^{-1}L(\theta).
\end{equation}
Let us end the proof of (a). We write it for $p\in(0,1]$; the case $p\geq 1$ is easier.

Suppose $\|x\|_p, \|y\|_p\leq 2^{-1/p}$, so that $\|x+y\|_p\leq 1$. Then 
$
 \| \Phi(x+y)- \Phi(x)- \Phi(y)\|_p^p$ is at most
\begin{gather*}
\| \Phi(x+y)- \Phi(x)- \Phi(y) -\Phi_0(x+y)+ \Phi_0(x)+ \Phi_0(y)\|_p^p+
  \| \Phi_0(x+y)- \Phi_0(x)- \Phi_0(y)\|_p^p\\
\leq  3\left(\frac{L(\theta)}{e}\right)^p+
\left(\frac{2L(\theta)}{e}(\||x|+|y|\|_p)
\right)^p = (3+2^p)\left(\frac{L(\theta)}{e}\right)^p,
\end{gather*}
so
$
\| \Phi(x+y)- \Phi(x)- \Phi(y)\|_p\leq 5^{1/p}e^{-1}L(\theta)$ provided $\|x\|_p, \|y\|_p\leq 2^{-1/p}).
$
But $\Phi$ is homogeneous and, therefore, for arbitrary $x,y\in\ell_p$ one has
$$
\| \Phi(x+y)- \Phi(x)- \Phi(y)\|_p\leq \frac{5^{1/p}}{e} L(\theta) 2^{1/p}\max(\|x\|_p,\|y\|_p)\leq \frac{10^{1/p}}{e}L(\theta)(\|x\|_p+\|y\|_p),
$$
as required.

Proof of (b), case $0<p\leq 1$. Suppose $\ell$ is a linear endomorphism of $\ell_p$ and set $\delta=\|\Phi-\ell\|$. Then since $\Phi(e_i)=0$ we have $\|\ell(e_i)\|_p\leq\delta$ and so
$$
\|\ell(s_n)\|_p^p=\|\ell(e_1)+\cdots+\ell(e_n)\|_p^p\leq n\delta^p.$$
But $\Phi(s_n)= \theta(\log n^{1/p})(s_n)$.
Combining,
$$
n(\theta(\log n^{1/p}))^p=\left\|\Phi(s_n)\right\|^p\leq \left\|\Phi(s_n)-\ell(s_n)\right\|^p+ n\delta^p\leq 2n\delta^p,
$$
hence $\delta\geq 2^{-1/p}\theta(\log n^{1/p})$ and, taking limit for $n\to\infty$ we see that the distance from $\Phi$ to the space of linear maps is at least $2^{-1/p}\|\theta\|_\infty$.

In order to get the corresponding estimate when $p\geq 1$ we observe that $\Phi$ maps $\ell_p^n$ to itself and since $\ell_p^n$ is complemented in $\ell_p$ by a norm-one projection, the distance from $\Phi$ to the space of linear endomorphisms of $\ell_p$ is at least the distance of the restriction of $\Phi$ to $\ell_p^n$ to the space of linear endomorphisms of the latter space.

The unit basis of $\ell_p^n$ is symmetric. This means that every operator of the form
\begin{equation*}
\sum_{i=1}^n\lambda_ie_i\longmapsto \sum_{i=1}^n\lambda_i\sigma_ie_{\pi(i)},
\end{equation*}
where $\pi$ is a permutation of $\{1,\dots,n\}$ and $\sigma_i=\pm1$ for $1\leq i\leq n$, is an isometry of $\ell_p^n$.

Let $U$ denote the group of all these operators. Then $\Phi$ commutes with $U$ in the sense that $\Phi\circ u=u\circ \Phi$ for all $u\in U$.  

On the other hand it is obvious that the only linear endomorphisms of $\ell_p^n$ that commute with $U$ are the scalar multiples of the identity.

Let $\ell$ be a linear map on $\ell_p^n$. Averaging  $\ell$ over $U$ we obtain the linear map ($|U|$ is the order of $U$)
$$
\tilde\ell(x)=\frac{1}{|U|}\sum_{u\in U}u^{-1}(\ell(u(x)))\quad\quad(x\in\ell_p^n).
$$
It turns out that $\|(\Phi-\tilde\ell):\ell_p^n\to\ell_p^n\|\leq \|\Phi-\ell\|$ since, given $x\in \ell_p^n$, one has 
\begin{align*}
\|\Phi(x)-\tilde\ell(x)\|_p
&=\left\| \frac{1}{|U|}\sum_{u\in U}u^{-1}(\Phi(u(x))) - \frac{1}{|U|}\sum_{u\in U}u^{-1}(\ell(u(x)))\right\|_p
\\
&\leq \frac{1}{|U|}\sum_{u\in U}\|u^{-1}(\Phi(u(x)))-u^{-1}(\ell(u(x)))\|\\
&= \frac{1}{|U|}\sum_{u\in U}\|\Phi(u(x))-\ell(u(x))\|\\
&\leq \frac{1}{|U|}\sum_{u\in U}\|\Phi-\ell\|\|u(x)\|_p\\
&= \|\Phi-\ell\|\|x\|_p.
\end{align*}
But $\tilde\ell$ commutes with $U$ and so it has to be a multiple of the identity. Hence it suffices to estimate the distances $\|\Phi-\alpha{\bf 1}_{\ell_p^n}\|$ for $\alpha\in \R$.
One has
$$
\|\Phi-\alpha{\bf 1}_{\ell_p^n}\|\geq \|\Phi(e_i)-\alpha e_i\|_p=|\alpha|.
$$
Also, taking again $s_n=e_1+\dots+e_n$, one has $\|s_n\|_p=n^{1/p}$ and $\Phi(s_n)=\theta(\log n^{1/p})s_n$, so
$$
\|\Phi-\alpha{\bf 1}_{\ell_p^n}\|\geq n^{-1/p}\|\Phi(s_n)-\alpha s_n\|_p=|\theta(\log n^{1/p})-\alpha|.
$$
Hence $\dist(\Phi, L(\ell_p^n))\geq {1\over 2}|\theta(\log n^{1/p})|$ and taking limit for $n\to\infty$ we see that
$$
\dist(\Phi,L(\ell_p))\geq \lim_{n\to\infty}\frac{\theta(\log n^{1/p})}{2}=\frac{\|\theta\|_\infty}{2}.
$$
This completes the proof.
\end{proof}

\subsection{Accessible operators arising by restriction or ``truncation''}\label{truncation}

We now present a simple procedure to get an accessible operator starting with a single nontrivial quasilinear map which works in many cases. We need the following concept, taken from \cite[p. 79]{k84}.

\begin{definition}
A quasi-Banach space is said to be an ultrasummand if it is complemented in all its ultrapowers through the diagonal embedding.
\end{definition}

A Banach space is an ultrasummand if and only if it is complemented in its bidual (or in any other dual space), see \cite[Proposition 1.4.13]{HMBST}. Thus for instance, $\ell_p$ and $L_p$ are ultrasummands for $1\leq p\leq \infty$, while $c_0$ is not. But there are nonlocally convex ultrasummands, among them $\ell_p$ and the Hardy classes $H_p$ for $0<p<1$. Neither $L_p$ not the $p$-Gurariy space $G_p$ are ultrasummands when $0<p<1$ as none of them is complemented in its ``countable'' ultrapowers. This last fact will not be used here and is far from being obvious: the point is that if $X$ denotes either $L_p$ or $G_p$ for a fixed $0<p<1$, then $X$ is separable, while there is no nonzero operator from $X_\U$ to any separable $p$-Banach space. See Notes 1.8.3 and 6.5.2 in \cite{HMBST}.

\medskip

Suppose $\phi:X\to Y$ is a nontrivial quasilinear map, with $Q[\phi]\leq 1$, and that $Y$ is an ultrasummand.
If $E$ is a finite-dimensional subspace of $X$, then $Q[\phi|_E]\leq 1$ and $\phi|_E$ is bounded. However,
$$
\sup\{\dist(\phi|_E,L(E,Y)): E \text{ is a finite-dimensional subspace of } X\}=\infty.
$$
For if not, there is a constant $D$ such that, for each finite-dimensional $E\subset X$ there is $\ell_E\in L(E, Y)$ satisfying $\|\phi|_E-\ell_E\|\leq D$.
Let $\mathscr F$ denote the set of all finite-dimensional subspaces of $X$, ordered by inclusion, and let $\mathscr V$ be an ultrafilter refining the order filter on $\mathscr F$ -- that is, $\mathscr V$ contains every set of the form $\{E: F\subset E\}$ for $F$ fixed in $\mathscr F$. Let us consider the mapping $L:X\to Y_\mathscr V$ given by $L(x)=(L_E(x))_\mathscr V$, where
$$
L_E(x)
=\begin{cases}\ell_E(x)&\text{if $x\in E$}\\
0& \text{otherwise}
\end{cases}
$$
It is easily seen that $L$ is a linear map: it is clearly homogeneous and, moreover, given $x,y\in X$, the set of those $E\in\mathscr F$ for which $L_E(x+y)=L_E(x)+L_E(y)$ belongs to $\mathscr V$ since it contains the ``tail'' $\{E\in\mathscr F: x,y\in E\}$.

Now, if $\Pi$ is a projection of $Y_\mathscr V$ onto $Y$, then for each $x\in X$, we have
\begin{align*}
\|\Pi(L(x))-\phi(x)\|&= \|\Pi(L(x))-\Pi\Delta\phi(x)\|\leq
\|\Pi\|\lim_{\mathscr V(E)}\|L_E(x)-\phi(x)\|\\&=
\|\Pi\|\lim_{\mathscr V(E)}\|\ell_E(x)-\phi(x)\|\leq
\|\Pi\|\cdot D\cdot\|x\|,
\end{align*}
that is, $\phi$ is trivial: a contradiction.

Therefore, we may take a sequence $(E_n)$ of finite-dimensional subspaces of $X$ such that $\dist(\phi|_{E_n}, L(E_n, Y))\to\infty$ as $n\to\infty$. For each $n$ choose $\ell_n\in  L(E_n, Y)$ such that $\|\phi|_{E_n}-\ell_n\|\leq \dist(\phi|_{E_n}, L(E_n, Y))+1/n$ and set
$$
\phi_n=\frac{\phi|_{E_n}-\ell_n}{\|\phi|_{E_n}-\ell_n\|}.
$$
\begin{example}
With the preceding notations, if $\U$ is any free ultrafilter on the integers, the family $(\phi_n)$ defines an accessible operator from $(E_n)_\U$ to $Y_\U$ that is not an ultraproduct of linear operators.
\end{example}

\begin{proof}
Since $Q[\phi_n]\leq 1/{\|\phi|_{E_n}-\ell_n\|}$ goes to zero, while $\|\phi_n\|=1$, the first part follows from Theorem~\ref{char}.
For the second we can estimate the distance between $\phi_n$ and the space of linear maps as follows:
\begin{align*}
\dist(\phi_n,L(E_n, Y))&= \frac{\dist(\phi|_{E_n}-\ell_n, L(E_n, Y))}{\|\phi|_{E_n}-\ell_n\|} =   \frac{\dist(\phi|_{E_n}, L(E_n, Y))}{\|\phi|_{E_n}-\ell_n\|}\\
&\geq \frac{\|\phi|_{E_n}-\ell_n\|-1/n}{\|\phi|_{E_n}-\ell_n\|}\geq 1-1/n,
\end{align*}
which completes the proof.
\end{proof}

\subsection{Accessible functionals on Banach spaces}\label{sec:AFBS}

Recall that a quasi-Banach space $X$ is said to be a $\K$-space if $K[X]$ is finite and that this happens if and only if every quasilinear functional $\phi:X\to \R$ is trivial, which in turn is equivalent to  ``every short exact sequence of quasi-Banach spaces of the form $0\to\R\to Z\to X\to 0$ splits''.

Which quasi-Banach spaces are $\K$-spaces? First of all: for $0<p\leq\infty$, the space $\ell_p$ (or $L_p$ over an arbitrary measure) is a $\K$-space if and only if $p\neq 1$.
Also, $B$-convex (Banach) spaces and $\mathscr{L}_\infty$-spaces are $\K$-spaces and so are quotients of Banach $\K$-spaces, see \cite{k78,kr}.

(A Banach space is said to be  $B$-convex if 
it has nontrivial type $p>0$; so $X$ is is not $B$-convex if and only if there is a constant $C$ such that for each $n\in\N$ there is an operator $u_n:\ell_1^n\to X$ such that $\|u_n\|\|u_n^{-1}\|\leq C$.) 

Here we want to exploit the last part of Corollary~\ref{K0}, namely that $X$ is a $\K_0$-space if and only if the accessible functionals on $X_\U$ are all in $X'_\U$.

Let us consider spaces of accessible functionals in more detail. Let $X_i$ be a family of Banach spaces, $\U$ a free ultrafilter on the index set $I$, and $(X_i)_\U$ the corresponding ultraproduct. Following Towsner, let us consider the following  subsets of the dual of $(X_i)_\U$:
\begin{itemize}
\item The ``trivial part'' $(X_i')_\U$.
\item The ``accessible part'' $\mathscr A((X_i)_\U)$.
\end{itemize}
After all these classes provided the initial motivation of this note.
Of course, one always has $(X_i')_\U \subset \mathscr A((X_i)_\U)\subset (X_i)_\U'$. 

In some cases, one has  $(X_i')_\U = (X_i)_\U'$: this happens if and only if $(X_i)_\U$ is reflexive or, equivalently, super-reflexive; see Heinrich \cite[Proposition 6.4]{heinrich}

As an immediate application, we have that super-reflexive Banach spaces are $\K_0$-spaces: indeed if $X$ is super-reflexive, then we have $(X_\U)'=(X')_\U$ and so every accessible linear functional is an ultraproduct of linear functionals. Actually the same idea was used in the one-sheet paper \cite{simple} (no joke intended)
to prove that super-reflexive Banach spaces are $\K$-spaces. We hasten to remark that super-reflexive Banach spaces are $B$-convex and so both results are contained in Kalton's for $B$-convex spaces.

Also, it is possible to have  $\mathscr A((X_i)_\U)=(X_i')_\U$ with $\mathscr A$ strictly smaller than $(X_i)_\U'$: this is the case if $X_i=X$, where $X$ is a nonreflexive Banach $\K_0$-space, for instance $X$ can be either $c_0$, a $C(K)$ space or any non-super-reflexive $B$-convex space.

We suspect that $\mathscr A((X_i)_\U)=((X_i)_\U)'$ forces $(X_i')_\U = (X_i)_\U'$, but we have a complete proof only for some particular cases.

\begin{lemma}
Let $\U$ be a free ultrafilter on $\N$. If $(X_n)$ is a sequence of finite-dimensional normed spaces such that $\mathscr A((X_n)_\U) = (X_n)_\U'$, then $(X_n')_\U = (X_n)_\U'$
\end{lemma}

\begin{proof}(Sketch)
If  $\mathscr A= (X_n')_\U$ there is nothing to prove. Otherwise the family $(X_n)$ cannot be ``uniformly $B$-convex'': there is a constant $C$, a sequence $k:\N\to\N$ such that $k(n)\to\infty$ along $\U$ and embeddings $u_n: \ell_1^{k(n)}\to X_n$ such that $\|u_n\|\,\|u_n^{-1}\|\leq C$. 

This implies that $(X_n)_\U$ contains a subspace isomorphic to $(\ell_1^{k(n)})_\U$, which in turn contains a complemented copy of $\ell_1(\langle k(n)\rangle_\U)$; see Lemma~\ref{lem:comple} below. Hence
$$
\dens ((X_n)_\U)' \geq \dens (\ell_1^{k(n)})_\U'\geq \dens \ell_1(\langle k(n)\rangle_\U)'=\dens \ell_\infty(2^{\aleph_0})=2^{2^{\aleph_0}}.
$$
On the other hand, since $X_n$ is finite-dimensional, the space of quasilinear functionals on $X_n$ is separable and
$$
\dens \mathscr A((X_n)_\U)\leq 2^{\aleph_0}< \dens (X_n)_\U',
$$
a contradiction.
\end{proof}

\begin{corollary}
Let $X$ be a Banach space with a basis and let $\U$ be a countably incomplete ultrafilter on $\N$. If each linear functional on $X_\U$ is accessible, then $X$ is super-reflexive and so $(X_\U)'=(X')_\U$.
\end{corollary}

\begin{proof}
For each $n\in\N$, let $X_n$ be the subspace spanned by the first $n$ vectors of the basis of $X$. The hypothesis implies that each linear functional on the ultraproduct $(X_n)_\U$ is accessible and the preceding lemma yields  $(X_n)_\U'=(X'_n)_\U$. Hence $(X_n)_\U$ is super-reflexive and so is its subspace $X$.
\end{proof}

\begin{remark}
There is an amusing interpretation of the accessible functionals in terms of internal subsets of the ultraproducts. Let $(X_i)_\U$ be an ultraproduct of Banach spaces. A subset $S$ of $(X_i)_\U$ is said to be  internal if there are subsets $S_i\subset X_i$ such that $S=(S_i)_\U$ in the sense that $x\in S$ if and only if it can be written as $(s_i)_\U$ with $s_i\in S_i$ for all $i\in I$. As it is well-known, each bounded linear functional corresponds to a closed hyperplane, the latter being the kernel of the former. A closed hyperplane of $(X_i)_\U$ corresponds to an accessible linear functional if and only if it is an internal subset. The point is that there are internal, closed hyperplanes that cannot be represented as an ultraproduct of  hyperplanes of the factors yet they can be represented as ultraproducts of ``hypersurfaces'' which are more and more ``flat''.
\end{remark}

\subsection{Accessible functionals on nonlocally convex spaces}\label{sec:AFNLCS}
We now present some examples of nonlocally convex $\K_0$-spaces. Here we take advantage of the fact that quasi-Banach spaces can have very few linear functionals. In fact they might have only the zero map: in this case we say that  the dual is trivial.

Let us say, for lack of better name, that a quasi-Banach space $X$ has continuous $L_p$-structure  if every $x\in X$ can be written as $x=x_1+x_2$, with $\|x\|^p=\|x_1\|^p+\|x_2\|^p$ and $\|x_1\|=\|x_2\|$ --- hence $\|x_i\|=2^{-1/p}\|x\|$ for $i=1,2$.

\begin{corollary}\label{Ex:K0}
Every quasi-Banach space whose countable ultrapowers have  trivial dual is a $\K_0$-space. In particular, for every $0<p<1$, the following are $\K_0$-spaces:
\begin{itemize}
\item[(a)] The space $L_p([0,1], E)$ for any quasi-Banach space $E$.
\item[(b)] The noncommutative $L^p$-space associated to a semifinite von Neumann algebra with no minimal projection.
\item[(c)] The $p$-Gurariy space $G_p$.
\end{itemize}
\end{corollary}

\begin{proof}
The first part trivially follows from the last part of Corollary~\ref{K0}.

As for the examples (a) and (b),
It is obvious that if $X$ has continuous $L_p$-structure then so do its ultrapowers and it is almost obvious that if $X$ has continuous $L_p$-structure for some $0<p<1$, then $X'=0$: indeed, if $f:X\to\R$ is linear and bounded then, given $x\in X$, write $x=x_1+x_2$ as in the definition, and observe that since $f(x)=f(x_1)+f(x_2)$ one has
$$
|f(x)|\leq 2^{1-1/p}\|f\|\|x\|;
$$
and this clearly implies that $\|f\|=0$.

Let us check that the spaces in (a) and (b) have continuous $L_p$-structure. For $L_p([0,1], E)$ this is trivial. Indeed, if $f:[0,1]\to E$ belongs to $L_p([0,1], E)$, then
$$\|f\|^p=\int_0^1 \|f(t)\|^p_E\: dt.$$
By the mean value theorem there is $s\in(0,1)$ such that 
$$\int_0^s \|f(t)\|^p_E\: dt= \int_s^1 \|f(t)\|^p_E\: dt $$
and $f_1=1_{[0,s]}f, f_2=1_{[s,1]}f$ does the trick.

\medskip

(b) We won't even define noncommutative $L^p$-spaces here. The reader is referred to the article by Pisier and Xu \cite{p-x} for details. The result follows from the fact that if $\mathscr M$ is semifinite and has no minimal projection, then every element of $L^p(\mathscr M)$ is contained in a subspace isometric to $L_p[0,1]$.
 (The hypothesis on $\mathscr M$ is necessary: if $\mathscr M$ is $B(H)$ ``itself'', then  
$L^p(\mathscr M)$ is the popular Schatten class $S_p$, which fails to be $\K$-spaces for $0<p<1$ by \cite[Theorem 6.1]{RMI}. Since $S_p$ has the BAP it cannot be a $\K_0$-space either.)

\medskip

(c) The $p$-Gurariy space $G_p$ is the only separable $p$-Banach space of ``almost universal disposition''. It was introduced in \cite{k78+} and is further studied in \cite{cgk} (see also \cite[Chapter~6]{HMBST}), where it is shown that its dual is zero and also that it is almost-transitive: given $x,y\in G_p$ such that $\|x\|=\|y\|=1$ and $\e>0$, there is a surjective linear isometry $u$ of $G_p$ such that $\|y-u(x)\|<\e$. Let us see that if $\U$ is a free ultrafilter on the integers, then $(G_p)_\U$ has trivial dual. Indeed, suppose $f:(G_p)_\U\to\R$ is a nonzero bounded linear functional and pick $(x_n)$ such that $f((x_n)_\U)\neq 0$. There is no loss of generality in assuming $\|x_n\|=1$ for all $n$. Fix a normalized $x\in G_p$ and, for each $n$ take a surjective linear isometry of $G_p$ such that $\|x_n-u_n(x)\|\leq 1/n$. The composition
$$
\begin{CD}
G_p@>\text{diagonal}>> (G_p)_\U@>(u_n)_\U>> (G_p)_\U@>f>>\R
\end{CD}
$$
is a nonzero linear functional on $G_p$, a contradiction.
\end{proof}

We continue with a couple of (related) examples. The first one shows that a nonzero linear functional can spring on the ultrapower of a quasi-Banach space with trivial dual. 

\begin{example}
A quasi-Banach space with trivial dual whose ultrapowers have nontrivial dual.
\end{example}

\begin{proof}
Set $X=\ell_1(\N,L_{p(n)})$, where $p(n)\to 1^-$ as $n\to\infty$ and $p(n)<1$ for all $n$. Then $X'=0$ since $\ell_1^0(\N,L_{p(n)})$ is dense in $X$ and $L_p$ has trivial dual for every $0<p<1$.

If $\U$ is a free ultrafilter on the integers, then $X_\U$ contains a complemented subspace isometric to $(L_{p(n)})_\U$, the ultraproduct of the family $(L_{p(n)})_n$ with respect to $\U$. But   $(L_{p(n)})_\U$ is a Banach space (actually it is linearly isometric to $L_1(\mu)$, for some very large measure $\mu$) since the quasinorm is subadditive: pick (uniformly) bounded families $f_n,g_n\in L_{p(n)}$. We have
\begin{align*}
\|(f_n)_\U+(g_n)_\U\|&=\lim_{\U(n)}\|f_n+g_n\|_{p(n)}\leq
\lim_{\U(n)}\left(\|f_n\|_{p(n)}^{p(n)}+\|g_n\|_{p(n)}^{p(n)}\right)^{1/{p(n)}}\\
&= \lim_{\U(n)}\left(\|f_n\|_{p(n)}^{p(n)}+\|g_n\|_{p(n)}^{p(n)}\right) =
\lim_{\U(n)}\|f_n\|_{p(n)}^{p(n)}+\lim_{\U(n)}\|g_n\|_{p(n)}^{p(n)}\\
&=\lim_{\U(n)}\|f_n\|_{p(n)}+\lim_{\U(n)}\|g_n\|_{p(n)}
=\|(f_n)_\U\|+\|(g_n)_\U\|.
\end{align*}
Hence $(L_{p(n)})_\U$ has nontrivial dual and so $X_\U$ does. 
\end{proof}

In the next example $X$ can be one of the spaces appearing in Corollary~\ref{Ex:K0}, but cannot be a Banach space!

\begin{example}
Let $X$ be a (nonzero) quasi-Banach space whose countable ultrapowers have trivial dual. Then the quotient of $X$ by any line (one-dimensional subspace) is a $\K_0$-space but not a $\K$-space.
\end{example}

\begin{proof}
Since ultrapowers preserve quotients and the property of having trivial dual is inherited by 
quotients the part concerning $\K_0$-spaces is obvious from Corollary~\ref{K0}. On the other hand, if $X'=0$, and $L\subset X$ is one-dimensional, the exact sequence
$$
\begin{CD}
0@>>> \R@>\imath>> X@>\pi >> X/L@>>> 0,
\end{CD}
$$
where $\imath$ takes $\R$ to $L$ and $\pi$ is the natural quotient map, does not split, since there is no linear operator from $X$ to $\R$, apart from zero. Hence $X/L$ is not a $\K$-space.
\end{proof}

We close this Section with a result on the accessible functionals on the ultrapowers of $\ell_p$ for $0<p<1$ which is roughly equivalent to a classical result by Kalton stating that those $\ell_p$ are $\K$-spaces; see \cite[Theorem 3.5]{k78}. The simple proof is based on the ``geometry'' of ultraproducts, and the key  point is
the following remark.

\begin{lemma}\label{lem:comple}
If $\U$ is a countably incomplete ultrafilter, then $(\ell_p)_\U$ contains a complemented subspace isometric to $\ell_p(\langle\N\rangle_\U)$ whose complement has continuous $L_p$-structure.
\end{lemma}

\begin{proof}
We first observe that, given $n:I\to \N$, the class of $(e_{n(i)})$ in $ (\ell_p)_\U$ depends only on the class of $(n(i))$ in the set-theoretic ultrapower $\langle\N\rangle_U$ and that this provides an isometric embedding of $\ell_p(\langle\N\rangle_U)$ into $(\ell_p)_\U$ by linear extension of the mapping
$$
J: e_{\langle n(i)\rangle_\U}\in \ell_p(\langle\N\rangle_U)\longmapsto (e_{n(i)})_\U\in (\ell_p)_\U.
$$
Let us check that the range of $J$ is a complemented subspace of $(\ell_p)_\U$. 
Given $x=(x_i)_\U$ in $(\ell_p)_\U$, we define a mapping $P(x):\langle\N\rangle_\U\to\R$ by the formula
$$
P(x)(\langle n(i)\rangle_\U)=\lim_{\U(i)}x_i(n(i)).
$$
It is easily seen that $\|P\|=1$ and that $P\circ J$ is the identity on $ \ell_p(\langle\N\rangle_\U)$. Moreover $x\in\ker P$ if and only if $\|x_i\|_\infty\to 0$ along $\U$. 
On the other hand, given $f\in \ell_p$ one always has a disjoint decomposition $f=f'+f''$, with 
$$
|\:\|f'\|^p- \|f''\|^p\:|\leq \|f\|_\infty^p.
$$
It follows that if $x=(x_i)_\U$ is normalized in $(\ell_p)_\U$ and $P(x)=0$, then we can take disjoint decompositions $x_i=x_i'+x_i''$ such that
$$
\lim_{\U(i)}\|x_i'\|_p=\lim_{\U(i)}\|x_i''\|_p=2^{-1/p}.
$$
Hence $\ker P$ has continuous $L_p$-structure.
\end{proof}

\begin{proposition}
Every accessible functional on $(\ell_p)_\U$ with $p\in(0,1)$ is an ultraproduct of linear functionals.
\end{proposition}

\begin{proof}
With the same notations as before, it is clear that every continuous linear functional on $(\ell_p)_\U$ vanishing on $\ell_p(\langle\N\rangle_\U)$ has to be zero: $\ker P$ has trivial dual.

Let $f_i:\ell_p\to \R$ be an admissible family, with accessible functional $f$. There is no loss of generality in assuming that each $f_i$ is homogeneous and that $\|f_i\|\leq \|f\|$ for every $i\in I$. 
For each $i\in I$, we define a linear functional on $\ell_p$ by the formula
$$
\ell_i(x)=\sum_n x(n)f_i(e_n).
$$
Note that $\|\ell_i\|\leq \|f_i\|$ and so the ultraproduct $(\ell_i)_\U$ is a linear functional on $(\ell_p)_\U$. But the difference  $(f_i)_\U-(\ell_i)_\U=(f_i-\ell_i)_\U$ is a linear functional vanishing on every $(e_{n(i)})_\U$ and so $f=(\ell_i)_\U$, as required.
\end{proof}

\subsection{The derivation out of space}
This section may be regarded as a small complement to Daws' \cite{daws}. We recommend the reader to have a look at Runde's monograph \cite{runde} to get an idea of amenability for Banach algebras and its many variations, connections and intricacies.
 Please note that our definitions are more akin to \cite{bade}, where weak amenability is introduced as a cheap substitute of amenability for commutative algebras.

If $p\in[1,\infty)$, the space $\ell_p$ is a commutative Banach algebra under the pointwise multiplication of functions, and the same occurs to $\ell_p^n$.

On the other hand,
if $A_i$ is a family of (commutative) Banach algebras indexed by $I$ and $\U$ is an 
 ultrafilter on $I$, then the Banach space ultraproduct $(A_i)_\U$ is a (commutative) Banach algebra under the coordinatewise product
$(a_i)_\U (b_i)_\U=(a_ib_i)_\U$. 
 This applies to $(\ell_p^n)_{\U}$ if $\U$ is a free ultrafilter on $\N$.

 Fix $p\in[1,\infty)$ and $n\in\N$, and consider the Kalton--Peck map $\Phi^n:\ell_p^n\to\ell_p^n$ given by $\Phi^n(x)=x\log(\|x\|_p/|x|)$, with the customary convention if $x(k)=0$ for some $1\leq k\leq n$. Kalton and Peck proved that $Q[\Phi^n]\leq 10^{1/p}e^{-1}$ --- this is ``our'' proof of the first part of Example~\ref{ex:kpmaps}. In particular $Q[\Phi^n]\leq 5$ for every $p\geq 1$ and every $n$. A direct application of Lagrange’s multiplier theorem shows that
 $$
 \sup\{ \|\Phi^n(x)\|_p : \|x\|_p\leq 1\} = \frac{\log n}{p}, 
 $$
 which is attained at $x=n^{-1/p}\sum_{k\leq n} e_k$, see \cite[Lemma 4]{ccs} for a more general result. It follows from Theorem~\ref{char} that if $\U$ is a free ultrafilter on $\N$ the map $D:(\ell_p^n)_{\U}\to (\ell_p^n)_{\U}$  defined by 
 \begin{equation}\label{eq:out}
D(x) = \left( \frac{x_n}{\log n}\cdot \log\frac{\|x_n\|_p}{|x_n|} \right)_\U
 \end{equation}
for $x=(x_n)_\U$ is an endomorphism of $(\ell_p^n)_{\U}$.  We want to see that $D$ is actually a derivation, with the meaning that it satisfies Leibniz rule $D(xy)=xDy+ yDx$. The reader should check that $D$ is not zero!

The estimate in \eqref{eq:estimate} shows that if $x_n\in\ell_p^n$ are such that $\sup_n\|x_n\|_p <\infty$, then 
$$
\left( \frac{x_n}{\log n}\cdot \log\frac{\|x_n\|_p}{|x_n|} \right)_n = \left( \frac{x_n}{\log n}\cdot \log\frac{1}{|x_n|} \right)_n
$$
in $(\ell_p^n)_\U$. 
Now, if $x=(x_n)_\U$ and $y=(y_n)_\U$, one has
\begin{align*}
D(xy)&=D\big{(}(x_ny_n)_\U\big{)}=\left( \frac{x_ny_n}{\log n} \log \frac{1}{|x_n y_n|}\right)_\U \\
&=
\left( x_n\frac{y_n}{\log n} \log\frac{1}{|y_n|}+ y _n\frac{x_n}{\log n} \log\frac{1}{|x_n|}\right)_\U = xDy + yDx. 
\end{align*}
This is a bit surprising because $\ell_p$ has no bounded derivation apart from zero. Actually $\ell_p$ is weakly amenable, meaning that if $X$ is a left-module over $\ell_p$ then every operator $D:\ell_p\to X$ satisfying Leibniz rule ($D(xy)=xDy+ yDx$) must be zero --- this is not the original definition, but an equivalent condition. And this is so because the idempotents (those elements satisfying $e^2=e$, think of the unit basis) generate a dense subspace of  $\ell_p$ and derivations annihilate idempotents: If $e^2=e$, then $De=2eDe$, so $eDe=2e^2De=2eDe=0$ and $De=0$. Note that $\ell_p$ fails to be amenable (in the usual sense) since it does not have an approximate unit.

We hasten to remark that the fact that $(\ell_p^n)_\U$ (and therefore the ultrapower algebra $(\ell_p)_\U$) fails to be weakly amenable follows from \cite[First part of the proof of Theorem 1.5]{bade} and has no particular interest
 in itself, it is the construction of the derivation on the ultraproduct what we find remarkable. 

\section{Concluding remarks}
We bid farewell to the (hypothesized) readers who have read this far with a slightly edited version of \cite[Note 5.4.1]{HMBST}.

The most interesting problems on Banach $\mathcal K$-spaces are to decide if every $\mathcal K_0$-space is a $\mathcal K$-space and to characterize Banach $\K$-
spaces by means of some easy to handle Banach space property: Kalton repeatedly conjectured that ``not containing $\ell_1^n$ uniformly complemented'' should be sufficient
; cf. \cite[p. 815]{kr}, \cite[p. 11]{k-mah}, \cite[Problem 4.2]{khandbook}. 
This would imply that any ultrapower (and so all even duals) of Banach $\K$-spaces are $\K$-spaces too and also that all Banach $\K_0$-spaces are $\K$-spaces. A first step in this direction appears in \cite[Corollary on p. 284]{HMBST}

Very little is known about the nature of Banach $\K$-spaces and, actually, the gap between Kalton's conjecture and the current list of members of the club of $\K$-spaces is sideral. In particular, we don't know whether the following are or are not $\K$-spaces:

\begin{itemize}
\item Pisier's spaces $P$ such that $P\otimes_\pi P= P\otimes_\varepsilon P$; cf. \cite{p}.
\item James' quasireflexive space in \cite{j50, j51}.
\item The original Tsirelson space $T^*$. 
\item The space of all (compact) bounded operators on a separable Hilbert space.
\item $L^p(\mathscr M)$ when $\mathscr M$ has no minimal projection and $0<p<1$.
\item The $p$-Gurariy space $G_p$ when $0<p<1$; see \cite[Section~7.5]{cgk}.
\item The Hardy classes $H_p$ for $0<p<1$; see
\cite[Problem 6]{k78}.
\item The disc algebra and the algebra of bounded analytic functions on the disc.
\item The spaces of vector-valued functions $\ell_p(E), L_p(E), c_0(E), C(K,E)$ when $p\neq 1$ and $E$ is a $\K$-space, as it is the case when $E=\ell_2$.
\item Put your favourite Banach space (containing $\ell_1^n$ uniformly, but not uniformly complemented) here. 
\end{itemize}

We know of no example of a Banach $\K$-space whose ultrapowers fail to be $\K$-spaces. And the same for $\K_0$-spaces and\,/\,or quasi-Banach spaces. On the other hand, if $X_\U$ is either a $\K_0$-space or a $\K$-space, then so is the base space $X$.
Thus, it is conceivable that  in the end  these classes are ``axiomatizable'' in the sense of  Henson--Iovino \cite[p. 82]{hen-iov}.

\section*{Acknowledgement}
I thank the anonymous referee (and so should the reader do) for the careful reading of the manuscript and suggestions that greatly improved the exposition.

\end{document}